\documentclass[11pt,b5paper,twoside,headrule]{amsart}
\usepackage[flushmargin]{footmisc}
\usepackage{amssymb}
\usepackage[all]{xy}
\newcommand{\di}{\displaystyle}
\newcommand{\intDelta}[1]{\ensuremath{\overset{\circ\hphantom{#1}}{\Delta^{#1}}}}

\newcommand{\st}{\mathcal{S}}
\newcommand{\A}{\mathbb{A}}
\newcommand{\B}{\mathbb{B}}
\newcommand{\LL}{\mathbb{L}}
\newcommand{\R}{\mathbb{R}}
\newcommand{\Z}{\mathbb{Z}}
\newcommand\Ab{\operatorname{Ab}}
\newcommand\Ob{\operatorname{Ob}}
\newcommand\Mod{\operatorname{Mod}}
\newcommand\Hom{\operatorname{Hom}}
\newcommand\Proj{\operatorname{Proj}}

\newtheorem{dummy}{realdumb}[section]

\newtheorem{lemma}[dummy]{Lemma}

\newtheorem{proposition}[dummy]{Proposition}
\theoremstyle{definition}           
\newtheorem{definition}[dummy]{Definition}
\newtheorem{example}[dummy]{Example}
\newtheorem{remark}[dummy]{Remark}
\theoremstyle{plain}

\begin{document}
\def\evenhead{{\protect\centerline{\textsl{\large{Andrew Ranicki and Michael Weiss}}}\hfill}}

\def\oddhead{{\protect\centerline{\textsl{\large{On The Algebraic $L$-theory of $\Delta$-sets}}}\hfill}}

\pagestyle{myheadings} \markboth{\evenhead}{\oddhead}

\thispagestyle{empty} 
\noindent{{\small\rm Pure and Applied
Mathematics Quarterly\\ Volume 8, Number 1\\ (\textit{Special
Issue: In honor of \\ F. Thomas Farrell and Lowell E. Jones})\\
1---29, 2012} \vspace*{1.5cm} \normalsize

\begin{center}
\Large{\bf On The Algebraic $L$-theory of $\Delta$-sets}
\end{center}

\begin{center}
{\large Andrew Ranicki and Michael Weiss}
\end{center}
\footnotetext{Received January 30, 2007.}


\bigskip

\begin{center}
\begin{minipage}{5in}
\noindent{\bf Abstract:} The algebraic $L$-groups $L_*(\A,X)$ are
defined for an additive category $\A$ with chain duality and a
$\Delta$-set $X$, and identified with the generalized homology
groups $H_*(X;\LL_{\bullet}(\A))$ of $X$ with coefficients in the
algebraic $L$-spectrum $\LL_{\bullet}(\A)$. Previously such groups
had only been defined for simplicial complexes $X$.
\end{minipage}
\end{center}

\section*{Introduction}

A `$\Delta$-set' $X$ in the sense of Rourke and Sanderson \cite{rs} is
a simplicial set without degeneracies. A simplicial complex
is a $\Delta$-set; conversely, the second barycentric (aka derived)
subdivision of a $\Delta$-set is a simplicial complex, and the homotopy
theory of $\Delta$-sets is the same as the homotopy theory of simplicial
complexes. However, $\Delta$-sets are sometimes more convenient than
simplicial complexes: they are generally smaller, and the quotient of
a $\Delta$-set by a group action is again a $\Delta$-set. In this paper
we extend the algebraic $L$-theory of simplicial complexes of Ranicki
\cite{ranicki1} to $\Delta$-sets.

In the original formulation of Wall \cite{wall} the surgery obstruction
theory of high-dimensional manifolds involved the algebraic $L$-groups
$L_*(R)$ of a ring with involution $R$, which are the Witt groups of
quadratic forms over $R$ and their automorphisms.  The subsequent
development of the theory in \cite{ranicki1} viewed $L_*(R)$ as the
cobordism groups of $R$-module chain complexes with quadratic
Poincar\'e duality, constructed a spectrum $\LL_{\bullet}(R)$ with
homotopy groups $L_*(R)$, and also introduced the algebraic $L$-groups
$L_*(R,X)$ of a simplicial complex $X$.  An element of $L_n(R,X)$ is a
cobordism class of directed systems over $X$ of $R$-module chain
complexes with an $n$-dimensional quadratic Verdier-type duality.  The
groups $L_*(R,X)$ were identified with the generalized homology groups
$H_*(X;\LL_{\bullet}(R))$, and the algebraic $L$-theory assembly map
$A:L_*(R,X) \to L_*(R[\pi_1(X)])$ was defined and extended to the
algebraic surgery exact sequence
$$\xymatrix@C-5pt{\dots \ar[r]& L_n(R,X) \ar[r]^-{\di{A}}&
L_n(R[\pi_1(X)]) \ar[r]& \st_n(R,X) \ar[r]& L_{n-1}(R,X) \ar[r]& \dots}$$
with $\st_n(R,X)$ the cobordism groups of the
$R[\pi_1(X)]$-contractible directed systems.  In particular, the
1-connective version gave an algebraic interpretation of the exact
sequence of the topological version of the
Browder-Novikov-Sullivan-Wall surgery theory: if the polyhedron $\Vert
X\Vert$ of a finite simplicial complex $X$ has the homotopy type of a
closed $n$-dimensional topological manifold then $\st_{n+1}(\Z,X)$ is
the structure set of closed $n$-dimensional topological manifolds $M$ with a
homotopy equivalence $M \simeq \Vert X\Vert$.

The Verdier-type duality of \cite{ranicki1} used the dual cells
in the barycentric subdivision of a simplicial complex $X$ to define
the dual of a directed system over $X$ of $R$-modules to be a
directed system over $X$ of $R$-module chain complexes.
The $\Delta$-set analogues of dual cells introduced by us in
Ranicki and Weiss \cite{ranickiweiss} are used here to define a
Verdier-type duality for directed systems of $R$-modules over a
$\Delta$-set $X$, which is used to define the generalized homology
groups $L_*(R,X)=H_*(X;\LL_{\bullet}(R))$ and an algebraic surgery
exact sequence as in the simplicial complex case.

The algebraic
$L$-theory of $\Delta$-sets is used in Macko
and Weiss \cite{mackoweiss}, and its multiplicative properties are
investigated in Laures and McClure \cite{lauresmcclure}.

\section{Functor categories\label{Functor categories}}
In this section, $X$ denotes a category with the following property. For every object $x$,
the set of morphisms to $x$ (with unspecified source) is finite; moreover, given morphisms $f:y\to x$
and $g:z\to x$ in $X$, there exists at most one morphism $h:y\to z$ such that $gh=f$.

Let $\A$ be an additive category with zero object $0 \in \Ob(\A)$.

\begin{definition}
(i) A function
$$M~:~\Ob(X) \to \Ob(\A)~;~x \mapsto M(x)$$
is {\it finite} if $M(x)=0$ for all but a finite number of objects $x$ in $\A$.\\
The direct sum
$\sum\limits_{x \in \Ob(X)}M(x)$ will be written as
$\sum\limits_{x \in X}M(x)$.\\
(ii) A functor $F:X \to \A$ is {\it finite} if the function
$F:\Ob(X)\to \Ob(\A)$ is finite.\hfill\qed
\end{definition}

\begin{definition}
(i)  The {\it contravariant functor category}
$\A_*[X]$ is the additive category of finite
contravariant functors $F:X \to \A$.
The morphisms in $\A_*[X]$ are the natural transformations.\\
(ii) The {\it covariant functor category} $\A^*[X]$ is
the additive category of covariant functors $F:X \to \A$.
The morphisms in $\A^*[X]$ are the natural transformations.
We write $\A^*_f[X]$ for the full subcategory whose objects are the
finite functors in $\A^*[X]$.
\hfill\qed
\end{definition}

\begin{remark}
We use the terminology $\A^*[X]$ for the
{\it covariant} functor category because it behaves contravariantly
in the variable $X$. Indeed a functor $g:X\to Y$ induces a functor
$\A^*[Y]\to \A^*[X]$ by composition with $g$. Our reasons for
using the terminology $\A_*[X]$ for the {\it contravariant}
functor category are similar, but more complicated. Below we
introduce a variation denoted $\A_*(X)$ which behaves
covariantly in $X$. \hfill\qed
\end{remark}

For the remainder of this section we shall only consider the contravariant
functor category $\A_*[X]$, but every result also
has a version for the covariant functor category $\A^*[X]$ (or $\A^*_f[X]$ in
some cases).

\begin{definition}
(i) A chain complex in an additive category $\A$
$$\xymatrix{C~:~\dots \ar[r] &
C_{n+1} \ar[r]^-{\di{d}} & C_n \ar[r]^-{\di{d}} & C_{n-1} \ar[r] &
\dots~~(d^2~=~0)}$$
is {\it finite} if $C_n=0$ for all but a finite number of $n \in \mathbb{Z}$.\\
(ii) Let $\B(\A)$ be the additive category of finite chain
complexes in $\A$ and chain maps.\hfill\qed
\end{definition}

A finite chain complex $C$ in $\A_*[X]$ is just an
object in $\B(\A)_*[X]$, and likewise for chain maps, so that
$$\B(\A_*[X])~=~\B(\A)_*[X]~.$$

\begin{definition} \label{weak}
A chain map $f:C \to D$ of chain complexes in $\A_*[X]$
is a {\it weak equivalence} if each
$$f[x]~:~C[x]~ \to~ D[x]~~(x \in X)$$
is a chain equivalence in $\A$.\hfill\qed
\end{definition}

A morphism $f:C \to D$ in $\B(\A_*[X])$ which is a chain equivalence is also a weak equivalence,
but in general a weak equivalence need not be a chain equivalence -- see
\ref{strong} for a more detailed discussion.

\begin{definition} \label{under}
Let $x$ be an object in $X$.\\
(i) The {\it under category} $x/X$ is the category with objects the
morphisms $f:x \to y$ in $X$, and morphisms $g:f \to f'$ the morphisms
$g:y \to y'$ in $X$ such that $gf=f'$
$$\xymatrix{ & x \ar[dl]_f \ar[dr]^{f'} \\ y \ar[rr]^g && y'}$$
The {\it open star} of $x$ is the set of objects in $x/X$
$$\text{st}(x)~=~\Ob(x/X)~=~\{ x \to y \}~.$$
(ii) The {\it over category} $X/x$ is the category with
morphisms $f:y\to x$ in $X$ as its objects, and so that morphisms $g:f \to f'$ are
the morphisms
$g:y \to y'$ in $X$ such that $f=f'g$
$$\xymatrix{y \ar[rr]^g \ar[dr]_f &  & y' \ar[dl]^{f'} \\ & x &}$$
The {\it closure} of $x$ is the set of objects in $X/x$
$$\text{cl}(x)~=~\Ob(X/x)~=~\{ y \to x \}~.$$
Because of our standing assumptions on $X$, the over category $X/x$
is isomorphic to a finite poset.
\hfill\qed
\end{definition}

In the applications of the contravariant functor category $\A_*[X]$
to topology we shall be particularly concerned with the subcategory of functors
satisfying the following property.

\begin{definition} A contravariant functor
$$F~:~X~ \to~ \A~;~x\mapsto F[x]$$
in $\A_*[X]$ is {\it induced} if there exists a finite function $x\mapsto F(x)\in \Ob(\A)$
and a natural isomorphism
\[  F[x] \cong \bigoplus_{x\to y}F(y)~. \]
The sum ranges over $\text{st}(x)$, and since the
function $x\mapsto F(x)$ is finite, $F[x]$ is only a sum of a finite
number of non-zero objects in $\A$. \newline
Similarly a covariant functor
$$F~:~X~ \to~ \A~;~x\mapsto F[x]$$
in $\A^*[X]$ is {\it induced} if there exists a function $x\mapsto F(x)\in \Ob(\A)$
and a natural isomorphism
\[  F[x] \cong \bigoplus_{y\to x}F(y)~. \]
\end{definition}

The full subcategories of the functor categories $\A_*[X]$, respectively $\A^*[X]$,
with objects the induced functors $F:X \to \A$
are equivalent, as we shall prove below, to the following categories.

\begin{definition}
Let $\A_*(X)$ be the additive category whose objects are
functions $x\mapsto F(x)$
such that $F(x)=0$ for all but a finite number of objects $x$.
A morphism $f:E \to F$ in $\A_*(X)$ is a collection of
morphisms $f(\phi):E(x) \to F(y)$ in $\A$, one for each
morphism $\phi:x \to y$ in $X$. The composite of the morphisms
$$f~=~\{f(\phi)\}~:~M\to N~~,~~g~=~\{g(\theta)\}~:~N\to P$$
is the morphism
$$gf~=~\{gf(\psi)\}~:~M\to P$$
with
$$gf(\psi:x \to z)~=~\sum\limits_{\phi:x \to y,\theta:y \to z,
\theta \phi=\psi}
g(\theta)f(\phi)~:~M(x)\to P(z)~.$$
We can view an object $F$ of $\A_*(X)$ as an object in $\A_*[X]$ by
writing
\[  F[x]=\bigoplus_{x\to y} F(y). \]
A morphism $\theta:w\to x$ in $X$ induces a morphism $F[x]\to F[w]$ in $\A$ which
maps the summand $F(y)$ corresponding to some $\phi:x\to y$ identically to the summand
$F(y)$ corresponding to the composition $\phi\theta:w\to y$. \newline
Let $\A^*(X)$ be the additive category whose objects are
functions $x\mapsto F(x)$. A morphism $f:E \to F$ in $\A_*(X)$ is a collection of
morphisms $f(\phi):E(y) \to F(x)$ in $\A$, one for each
morphism $\phi:x \to y$ in $X$. Again we can view an object $F$ of $\A^*(X)$ as an object in $\A^*[X]$ by
writing
\[  F[x]=\bigoplus_{y\to x} F(y). \]
\end{definition}

\begin{proposition} \label{adjoint}
{\rm (i)} For any object $M$ in $\A_*(X)$ and any object $N$
in $\A_*[X]$
$$\Hom_{\A_*[X]}(M,N)~=~
\sum\limits_{x \in X}\Hom_{\A}(M(x),N[x])~.$$
{\rm (ii)} For any objects $L,M$ in $\A_*(X)$
$$
\Hom_{\A_*[X]}(L,M)~=~\sum\limits_{x \to y}
\Hom_{\A}(L(x),M(y))~.$$
{\rm (iii)} The additive category $\A_*(X)$ is equivalent
to the full subcategory of the contravariant functor category
$\A_*[X]$ with objects the induced functors.
\end{proposition}

\begin{proof}
(i) A morphism $f:M \to N$ in  $\A_*[X]$ is determined by the
composite morphisms in $\A$
$$\xymatrix@C+15pt{M(x) \ar[r]^-{\di{\text{inclusion}}} &
M[x] \ar[r]^-{\di{f[x]}} & N[x]~~(x \in X)}$$
(ii) By (i), a morphism $f:L \to M$ in  $\A_*[X]$ is determined by the
composite morphisms in $\A$
$$\xymatrix@C+15pt{L(x) \ar[r]^-{\di{\text{inclusion}}} &
L[x] \ar[r]^-{\di{f[x]}} & \rule{0mm}{6.5mm} M[x]~=~\sum\limits_{x\to y}M(y)
~~(x \in X)~.}$$
(iii) Every object $M$ in $\A_*(X)$
determines an induced contravariant functor
$$X\to \A~;~x\mapsto M[x]~=~\sum\limits_{x \to y}M(y)~,$$
i.e. an object in $\A_*[X]$,
and every induced functor is naturally equivalent to one of this type.
\end{proof}

\begin{proposition} \label{local}
The following conditions on a chain map $f:C \to D$
in $\A_*(X)$ are equivalent:
\begin{itemize}
\item[\rm (a)] $f$ is a chain equivalence,
\item[\rm (b)]  each of the component chain maps in $\A$
$$f(1_x)~:~C(x)\to  D(x)~~(x \in X)$$
is a chain equivalence,
\item[\rm (c)]  $f:C\to D$ is a weak equivalence in $\A_*[X]$, that is,
$C[x]\to D[x]$ is a chain equivalence for all $x$.
\end{itemize}
\end{proposition}
\begin{proof} The proof given in Proposition 2.7 of Ranicki and Weiss \cite{ranickiweiss}
in the case when $\A$ is the additive category of $R$-modules (for
some ring $R$) works for an arbitrary additive category. 
\end{proof}

\begin{remark} \label{strong}
Every chain equivalence of chain complexes in $\A_*[X]$
is a weak equivalence. By \ref{local} every weak equivalence of degreewise induced
finite chain complexes in $\A_*[X]$ is a chain equivalence.
See Ranicki and Weiss \cite[1.13]{ranickiweiss} for an explicit example of a
weak equivalence of finite chain complexes in $\A_*[X]$
which is not a chain equivalence. It is proved in \cite[2.9]{ranickiweiss} that
every finite chain complex $C$ in $\A_*[X]$ is weakly equivalent
to one in $\A_*(X)$. \hfill\qed
\end{remark}

\section{$\Delta$-sets\label{Delta-sets}}

Let $\Delta$ be the category with objects the sets
$$[n]~=~\{0,1,\dots,n\}~~(n \geqslant 0)$$
and morphisms $[m] \to [n]$ order-preserving injections.
Every such morphism has a unique factorization as the composite of
the order-preserving injections
\begin{center}
$\partial_i~:~[k-1]\to [k]~;~j \mapsto
\begin{cases}
j &\text{if $j <i$} \cr
j+1 &\text{if $j\geqslant i$}~.
\end{cases}$
\end{center}
\begin{definition} (Rourke and Sanderson \cite{rs})
A {\it $\Delta$-set} is a contravariant functor
$$X~:~\Delta\to \{\text{sets and functions}\}~;~[n]\mapsto X^{(n)}~.$$
\hfill\qed
\end{definition}

Equivalently, a $\Delta$-set $X$ can be regarded as
a sequence $X^{(n)}$ $(n\geqslant 0)$ of sets, together with face maps
$$\partial_i~:~ X^{(n)}~\to ~ X^{(n-1)}~~  (0\leqslant i\leqslant n)$$
such that
$$\partial_i\partial_j~ =~ \partial_{j-1}\partial_i~~\hbox{for } i<j~ .$$
The elements $x \in X^{(n)}$ are the {\it $n$-simplices} of $X$.

\begin{definition} (Rourke and Sanderson \cite{rs})\\
(i) The {\it realization} of a $\Delta$-set $X$ is the $CW$ complex
$$\Vert X \Vert~=~\coprod\limits_{n=0}^{\infty}(X^{(n)} \times \Delta^n)/\sim$$
with
$$\begin{array}{l}
\Delta^n~=~\{(s_0,s_1,\dots,s_n)\in \R^n\,\vert\,0 \leqslant s_i \leqslant 1,~
\sum\limits^n_{i=0}s_i=1\}~,\\
\partial_i~:~\Delta^{n-1} \hookrightarrow \Delta^n~;~(s_0,s_1,\dots,s_{n-1})
\mapsto (s_0,s_1,\dots,s_{i-1},0,s_{i+1},\dots,s_n)~,\\[1ex]
(x,\partial_is) \sim (\partial_i x,s)~~(x \in X^{(n)},s \in \vert\Delta^{n-1}\vert)~.
\end{array}$$
(ii) There is one $n$-cell $x(\Delta^n) \subseteq \Vert X \Vert$
for each $n$-simplex $x \in X$, with {\it characteristic map}
$$x~:~\Delta^n \to \Vert X \Vert~;~(s_0,s_1,\dots,s_n)
\mapsto (x,(s_0,s_1,\dots,s_n))~.$$
The {\it boundary}  $x(\partial \Delta^n) \subseteq \Vert X \Vert$ is
the image of
$$\begin{array}{ll}
\partial \Delta^n&=~\bigcup\limits^n_{i=0} \partial _i\vert \Delta^{n-1}\vert\\
&=~\{(s_0,s_1,\dots,s_n)\in \R^n\,\vert\,
0 \leqslant s_i \leqslant 1,~
\sum\limits^n_{i=0}s_i=1,~\hbox{\rm $s_i=0$ for some $i$}\}
\end{array}$$
and the {\it interior} $x(\intDelta{n}) \subseteq \Vert X \Vert$
is the image of
$$\begin{array}{ll}
\intDelta{n}&=~
\Delta^n\backslash \partial\Delta^n\\[1ex]
&=~\{(s_0,s_1,\dots,s_n)\in \R^n\,\vert\,0 < s_i \leqslant 1,~
\sum\limits^n_{i=0}s_i=1\} \subseteq  \Delta^n ~.
\end{array}$$
The characteristic map $x: \Delta^n \to \Vert X \Vert$
is injective on $ \intDelta{n}\subseteq  \Delta^n$.\hfill\qed
\end{definition}
\begin{example}  Let $\Delta^n$ be the $\Delta$-set with
$$(\Delta^n)^{(m)}~=~\{\hbox{morphisms}~[m] \to [n]~\hbox{in}~\Delta\}
~~(0 \leqslant m \leqslant n)~.$$
The realization $\Vert \Delta^n\Vert$ is the geometric $n$-simplex $ \Delta^n $
(as in the above definition). It should be clear from the context whether
$\Delta^n$ refers to the $\Delta$-set or the geometric realization.\hfill\qed
\end{example}

We regard a $\Delta$-set $X$ as a category, whose objects are the simplices,
writing the dimension of an object $x \in X$ as $\vert x \vert$, i.e.
$\vert x \vert =m$ for $x \in X^{(m)}$.
A morphism $f:x \to y$ from an $m$-simplex $x$ to an $n$-simplex $y$
is a morphism $f:[m] \to [n]$ in $\Delta$ such that
$$f^*(y)~=~x \in X^{(m)}~.$$
In particular, for any $x \in X^{(m)}$ with $m \geqslant 1$
there are defined $m+1$ distinct morphisms in $X$
$$\partial_i~:~\partial_ix \to x~~(0 \leqslant i \leqslant m)~.$$

\begin{example}\label{expl-roundADelta} (i) Let $X$ be a $\Delta$-set. An object $M$ of $\A_*(X)$ is just an object
$M$ of $\A$ with a direct sum decomposition $M=\bigoplus_{x\in X}M(x)$. A morphism $f:M\to N$
in $\A_*(X)$ is a collection of morphisms $f_{xy,\lambda}: M(x)\to N(y)$, one such for every pair of
simplices $x,y$ and face operator $\lambda$ such that $\lambda^*y=x$. \newline
We like to think of a morphism $f:M\to N$ in $\A_*(X)$ as a morphism in $\A$ with additional structure.
Source and target of that morphism in $\A$ are $M(X)=\bigoplus_x M(x)$ and $N(X)=\bigoplus_x N(x)$,
respectively. For simplices $x$ and $y$, the $xy$-component of the morphism $M(X)\to N(X)$ determined by $f$ is
\[  \sum_{\lambda} f_{xy,\lambda} \]
where the sum runs over all $\lambda$ such
 that $\lambda^*y=x$. \newline
(ii) If $X$ is a simplicial complex then a morphism
in $\A_*(X)$ is just a morphism $f:M \to N$ in $\A$ between objects with
finite direct sum decompositions
$$M~=~\sum\limits_{x \in X}M(x)~,~N~=~\sum\limits_{y \in X}N(y)$$
such that the components $f(x,y):M(x) \to N(y)$ are 0 unless $x \leqslant y$. \newline
(iii) The description of $\A_*(X)$ in (ii) also applies in the case of a $\Delta$-set $X$
where, for any two simplices $x$ and $y$, there is at most one morphism from $x$ to $y$.
In particular it applies when $X=Y'$ is the barycentric subdivision of another $\Delta$-set $Y$,
to be defined in the next section.
\hfill\qed
\end{example}

\begin{definition}\label{chain}
Let $X$ be a $\Delta$-set, and let $R$ be a ring.\\
(i) The {\it $R$-coefficient simplicial chain complex of $X$} is the
free (left) $R$-module chain complex $\Delta(X;R)$ with
$$d~=~\sum\limits^n_{i=0}(-)^i\partial_i~:~
\Delta(X;R)_n~=~R[X^{(n)}]\to \Delta(X;R)_{n-1}~=~R[X^{(n-1)}]~.$$
The {\it $R$-coefficient homology} of $X$ is the homology of $\Delta(X;R)$
$$H_*(X;R)~=~H_*(\Delta(X;R))~=~H_*(\Vert X \Vert;R)~,$$
noting that $\Delta(X;R)$ is the $R$-coefficient cellular chain complex of
$\Vert X \Vert$.\\
(ii) Suppose that $R$ is equipped with an involution
$$R \to R~;~r \mapsto \overline{r}$$
(e.g. the identity for a commutative ring), allowing the definition
of the {\it dual} of an $R$-module $M$ to be the $R$-module
$$M^*~=~\Mod_R(M,R)~,~
R \times M^* \to M^*~;~(r,f) \mapsto (x \mapsto f(x) \overline{r})~.$$
The {\it $R$-coefficient simplicial cochain complex of $X$}
$$\Delta(X;R)^*~=~{\rm Hom}_R(\Delta(X;R),R)$$
is the $R$-module cochain complex with
$$d^*~=~\sum\limits^{n+1}_{i=0}(-)^i\partial^*_i~:~
\Delta(X;R)^n~=~R[X^{(n)}]^*\to \Delta(X;R)^{n+1}~=~R[X^{(n+1)}]^*~,$$
The {\it $R$-coefficient cohomology} of $X$ is the cohomology of
$\Delta(X;R)^*$
$$H^*(X;R)~=~H^*(\Delta(X;R)^*)~=~H^*(\Vert X \Vert;R)~,$$
noting that $\Delta(X;R)^*$ is the $R$-coefficient cellular cochain complex of
$\Vert X \Vert$.\hfill\qed
\end{definition}

A simplicial complex $X$ is {\it ordered} if the vertices
in any simplex are ordered, with faces having compatible orderings.
From now on, in dealing with simplicial complexes we shall
always assume an ordering.

\begin{example}
A simplicial complex $X$ can be regarded as a
$\Delta$-set, with $X^{(n)}$ the set of $n$-simplices and
$$\partial_i~:~X^{(n)}\to X^{(n-1)}~;~(v_0v_1\dots v_n)\mapsto
(v_0v_1 \dots v_{i-1} v_{i+1} \dots v_n)~.$$
There is one morphism $x \to y$ in $X$ for each face inclusion $x \leqslant y$.
The realization $\Vert X \Vert$ of $X$ regarded as a $\Delta$-set
is the polyhedron of the simplicial complex $X$, with the characteristic maps
$x:\Delta^{\vert x \vert} \to \Vert X \Vert$ ($x \in X$) injections.
The simplicial chain complex $\Delta(X;R)$ is just the usual $R$-coefficient
simplicial chain complex of $X$, and $\Delta(X;R)^*$ is the
$R$-coefficient simplicial cochain  complex of $X$. \hfill\qed
\end{example}
\begin{example} \label{over}
Let $X$ be a $\Delta$-set, and let $x \in X$ be a simplex.\\
(i) In general, the canonical map
$$\Ob(x/X)~=~\text{st}(x) \to \Ob(X)~;~(x \to y) \mapsto y$$
is not injective.
The simplices $y \in \Ob(X)\backslash \text{im}(\text{st}(x))$ are the
objects of a sub-$\Delta$-set $X\backslash \text{im}(\text{st}(x)) \subset X$.
If $X$ is a simplicial complex then
$\text{st}(x) \to \Ob(X)$ is injective,
and $X\backslash \text{st}(x) \subset X$ is the subcomplex
with simplices $y \in X$ such that $x \not\leqslant y$.\\
(ii) The over category $X/x=\{y \to x\}$ (\ref{under}) is a $\Delta$-set with
$$(X/x)^{(n)}~=~\{y \to x\,\vert\, y \in X^{(n)}\}~~(n \geqslant 0)~.$$
It is isomorphic as a $\Delta$-set to $\Delta^{|x|}$.
The forgetful functor
$$X/x \to X~;~(y \to x) \mapsto y$$
is a $\Delta$-map, inducing the characteristic map $\Delta^{|x|}\to \Vert X \Vert$.
If $X$ is a simplicial complex then $X/x \to X$ is injective, and so is the
induced characteristic map. \hfill\qed
\end{example}

\begin{example} \label{circle1}
(i) If a group $G$ acts on a $\Delta$-set $X$ the quotient
$X/G$ is again a $\Delta$-set, with realization $\Vert X/G\Vert=\Vert X \Vert/G$.
However, if $X$ is a simplicial complex and $G$ acts on $X$, then
$X/G$ is not in general a simplicial complex. See (ii) for an
example.\\
(ii) Suppose $X=\R$, the $\Delta$-set with
$$X^{(0)}~=~X^{(1)}~=~\Z~,~\partial_0(n)=n~,~\partial_1(n)=n+1~,$$
and let the infinite cyclic group $G=\Z=\{t\}$ act on  $X$ by $tn=n+1$.
The quotient $\Delta$-set $S^1=\R/\Z$ is the circle, with one
0-simplex $x_0$ and one 1-simplex $x_1$
$$(S^1)^{(0)}~=~\{x_0\}~,~(S^1)^{(1)}~=~\{x_1\}~,~
\partial_0(x_1)~=~\partial_1(x_1)~=~x_0~.$$
\hfill\qed
\end{example}

\begin{example} For any space $M$ use the standard $n$-simplices
$\Delta^n$ and face inclusions $\partial_i:\Delta^{n-1} \hookrightarrow \Delta^n$
to define the {\it singular $\Delta$-set} $X=M^{\Delta}$ by
$$X^{(n)}~=~M^{\Delta^n}~,~\partial_i~:~X^{(n)} \to X^{(n-1)}~;~x \mapsto
x \circ \partial_i~.$$
We shall say that a singular simplex $x:\Delta^n \to X$
is a face of a singular simplex $y:\Delta^m \to X$ if
$x=y\circ \partial_{i_1}\circ\dots\circ\partial_{i_{m-n}}$ for a given face inclusion
$$\partial_{i_1}\circ\dots\circ\partial_{i_{m-n}}~:~\Delta^n \hookrightarrow \Delta^m~,$$
writing $x \leqslant y$ (and $x<y$ if $x \neq y$).
The simplicial chain complex $\Delta(X;R)=S(M;R)$ is just the usual
$R$-coefficient singular chain complex of $M$, so that
$$H_*(\Vert X \Vert;R)~=~H_*(X;R)~=~H_*(M;R)~.$$
Also $\Delta(X;R)^*=S(M;R)^*$ is the $R$-coefficient
singular cochain  complex of $M$, and
$$H^*(\Vert X \Vert;R)~=~H^*(X;R)~=~H^*(M;R)~.$$
\hfill\qed
\end{example}

\section{The barycentric subdivision \label{barycentric}}

The $\Delta$-set analogue of the barycentric
subdivision $X'$ of a simplicial complex $X$ and the dual cells
$D(x,X) \subset X'$ ($x \in X$)
makes use of the following standard categorical construction.

\begin{definition}
(i) The {\it nerve} of a category $\mathcal{C}$ is the
simplicial set with one $n$-simplex for
each string $x_0 \to x_1 \to \dots \to x_n$ of morphisms in $\mathcal{C}$,
with
$$\partial_i(x_0 \to x_1 \to \dots \to x_n)~=~
(x_0 \to x_1 \to \dots \to x_{i-1} \to x_{i+1} \to \dots \to x_n)~.$$
(ii) An $n$-simplex $x_0 \to x_1 \to \dots \to x_n$ in the nerve is
{\it non-degenerate} if none of the morphisms $x_i \to x_{i+1}$ is
the identity.\hfill\qed
\end{definition}

If the category $\mathcal{C}$ has the property that the composite of
non-identity morphisms is a non-identity, then the non-degenerate
simplices in the nerve define a $\Delta$-set, which we shall also call
the nerve and denote by $\mathcal{C}$.

\begin{definition}\label{dual}
(Rourke and Sanderson \cite[\S4]{rs}, Ranicki and Weiss \cite[1.6, 1.7]{ranickiweiss})
Let $X$ be a $\Delta$-set.\\
(i) The {\it barycentric subdivision} of $X$ is the $\Delta$-set
$X'$ defined by the nerve of the category $X$.\\
(ii) The {\it dual} $x^{\perp}$ of a simplex $x\in X$ is the
nerve of the under category $x/X$ (\ref{under}).
An $n$-simplex in the $\Delta$-set $x^{\perp}$ is thus a sequence
of morphisms in $X$
$$x \to x_0 \to x_1 \to \dots \to x_n$$
such that $x_0 \to x_1 \to \dots \to x_n$ is non-degenerate. In particular
$$(x^{\perp})^{(0)}~=~\{x \to x_0\}~=~\text{st}(x)~.$$
(iii) The {\it boundary of the dual} $\partial x^{\perp}$ is the
sub-$\Delta$-set of $x^{\perp}$ consisting of the $n$-simplices
$x \to x_0 \to x_1 \to \dots \to x_n$
such that $x \to x_0$ is not the identity.\hfill\qed
\end{definition}

The under category $x/X$ has an initial object, so that the nerve $x^{\perp}$ is
contractible. The rule $x \to x^{\perp}$ is contravariant, i.e. every
morphism $x \to y$ induces a $\Delta$-map $y^{\perp} \to x^{\perp}$.

\begin{lemma} The realizations $\Vert X \Vert$, $\Vert X' \Vert$ of a
$\Delta$-set $X$ and its barycentric subdivision $X'$ are homeomorphic,
via a homeomorphism $\Vert X' \Vert \to \Vert X \Vert$ sending the
vertex $x \in X=(X^{\prime})^{(0)}$ to the barycentre
$$\widehat{x}~=~x(\dfrac{1}{n+1},\dfrac{1}{n+1},\dots,\dfrac{1}{n+1})
\in x(\intDelta{n}) \subseteq \Vert X \Vert~.$$
\end{lemma}
\begin{proof} It suffices to consider the special case $X=\Delta^n$,
so that $X$ and $X'$ are simplicial complexes,
and to define a homeomorphism
$\Vert X' \Vert \to \Vert X \Vert$ by $x \mapsto \widehat{x}$
and extending linearly.
\end{proof}
\begin{definition} Let $X$ be a $\Delta$-set, and let $x \in X$ be a
simplex.\\
(i)  The {\it open star space}
$$\Vert\text{st}(x)\Vert~=~\bigcup\limits_{y \in x^{\perp}\backslash\partial x^{\perp}}
\intDelta{\vert y \vert}\subseteq \Vert X' \Vert~=~\Vert X \Vert$$
is the subspace of the realization $\Vert X' \Vert$ of the barycentric
subdivision $X'$ defined by the union of the interiors of the simplices
$y \in x^{\perp}\backslash\partial x^{\perp}$, i.e.
$$y~=~(x\to x_0 \to \dots \to x_n) \in X'$$
with $x \to x_0=x$ the identity.\\
(ii) The {\it homology} of the open star is
$$H_*(\text{st}(x))~=~H_*(\Delta(\text{st}(x)))$$
with $\Delta(\text{st}(x))$ the chain complex defined by
$$\Delta(\text{st}(x))~=~\Delta(x^{\perp},\partial x^{\perp})_{*-\vert x\vert}~.$$
\hfill\qed
\end{definition}

\begin{lemma}
For any simplex $x \in X$ of a $\Delta$-set $X$ the
characteristic $\Delta$-map
$$i~:~x^{\perp} \to X'~;~(x \to x_0 \to \dots \to x_n) \mapsto
(x_0 \to \dots \to x_n)$$
is injective on $x^{\perp}\backslash \partial x^{\perp}$.
The images $i(\partial x^{\perp}),i(x^{\perp})\subseteq X'$ are
sub-$\Delta$-sets such that
$$\Vert i(x^{\perp})\Vert \backslash
\Vert i(\partial x^{\perp})\Vert~=~\Vert \hbox{\rm st}(x) \Vert  \subseteq
\Vert X' \Vert$$
and there are homology isomorphisms
$$\begin{array}{ll}
H_*({\rm st}(x))&=~H_{*-\vert x \vert}(x^{\perp},\partial x^{\perp})\\[1ex]
&\cong~H_{*-\vert x \vert}(i(x^{\perp}),i(\partial x^{\perp}))\\[1ex]
&\cong~H_*(\Vert X \Vert,\Vert X \Vert \backslash \Vert \hbox{\rm st}(x)
\Vert)\\[1ex]
&\cong~
H_*(\Vert X \Vert,\Vert X \Vert\backslash \{\widehat{x}\})~.
\end{array}$$
\end{lemma}
\begin{proof} The inclusion
$(\Vert X \Vert,\Vert X \Vert\backslash \Vert \hbox{\rm st}(x)\Vert) \hookrightarrow
(\Vert X \Vert, \Vert X \Vert \backslash \{\widehat{x}\})$
is a deformation retraction, and the open star subspace
$\Vert\hbox{\rm st}(x)\Vert \subset \Vert X \Vert$
has an open regular neighbourhood
\[\Vert\hbox{\rm st}(x)\Vert \times \intDelta{\vert x \vert}\subset \Vert X \Vert\]
with one-point compactification
$$(\Vert\hbox{\rm st}(x)\Vert \times \intDelta{\vert x \vert})^{\infty}~=~
\Vert i(x^{\perp})\Vert/\Vert i(\partial x^{\perp})\Vert \wedge \Delta^{\vert x \vert}/\partial \Delta^{\vert x \vert}~,$$
so that
$$\begin{array}{ll}
H_*(\Vert X \Vert,\Vert X \Vert\backslash \{\widehat{x}\})&\cong~
H_*(\Vert X \Vert,\Vert X \Vert \backslash \Vert \hbox{\rm st}(x) \Vert)\\[1ex]
&\cong~\widetilde{H}_*(\Vert i(x^{\perp})\Vert/\Vert
i(\partial x^{\perp})\Vert \wedge \Delta^{\vert x \vert}/\partial \Delta^{\vert x \vert})\\[1ex]
&\cong~H_{*-\vert x \vert}(i(x^{\perp}),i(\partial x^{\perp}))~.
\end{array}$$
\end{proof}

\begin{example} \label{dual cell}
Let $X$ be a simplicial complex. The barycentric subdivision of $X$ is the
ordered simplicial complex $X'$ with one $n$-simplex for each sequence of proper face
inclusions $x_0 < x_1 < \dots < x_n$. By definition,
the {\it dual cell} of a simplex $x \in X$
is the subcomplex $D(x,X) \subseteq X'$ consisting of all the simplices
$x_0 < x_1 < \dots < x_n$ with $x \leqslant x_0$. The {\it boundary} of the dual
cell is the subcomplex $\partial D(x,X) \subseteq D(x,X)$
consisting of all the simplices $x_0 < x_1 < \dots < x_n$ with $x < x_0$.
The $\Delta$-sets associated to $X,X',D(x,X),\partial D(x,X)$ are just
the $\Delta$-sets $X,X',x^{\perp},\partial x^{\perp}$ of \ref{dual},
with the characteristic map $i:x^\perp=D(x,X) \to X'$ injective.
Moreover, $X\backslash {\rm st}(x) \subset X$ is a subcomplex such that
$$\Vert X\backslash {\rm st}(x) \Vert ~=~\Vert X \Vert \backslash \Vert {\rm st}(x)\Vert$$
and
$$\Delta({\rm st}(x))~=~\Delta(D(x,X),\partial D(x,X))_{*-\vert x \vert}~\simeq~
\Delta(X,X \backslash {\rm st}(x))~.$$
\hfill\qed
\end{example}

\begin{example} \label{circle2}
Let $X$ be the $\Delta$-set (\ref{circle1})
with one 0-simplex $x_0$ and one 1-simplex $x_1$,
with non-identity morphisms
$$\xymatrix@C-5pt{x_0 \ar@<0.75ex>[r]  \ar@<-0.75ex>[r] & x_1}$$
and realization $\Vert X \Vert=S^1$.
The barycentric subdivision $X'$ is the $\Delta$-set with 2 0-simplices and 2
1-simplices:
$${X'}^{(0)}~=~\{x_0,x_1\}~, ~
{X'}^{(1)}~=~\{\xymatrix@C-5pt{x_0 \ar@<0.75ex>[r]  \ar@<-0.75ex>[r] & x_1}\}~.$$
The duals and their boundaries are given by
$$\begin{array}{l}
x_0^{\perp}~=~\{\xymatrix@C-5pt{x_0\ar[r] &x_0},~
\xymatrix@C-5pt{x_0 \ar@<0.75ex>[r]  \ar@<-0.75ex>[r] & x_1}\}
\cup\{\xymatrix@C-5pt{x_0 \ar[r] &x_0 \ar@<0.75ex>[r]  \ar@<-0.75ex>[r] & x_1}\}~,\\[2ex]
\partial x_0^{\perp}~=~\{
\xymatrix@C-5pt{x_0 \ar@<0.75ex>[r]  \ar@<-0.75ex>[r] & x_1}\}~=~\{0,1\}~,\\[2ex]
x_1^{\perp}~=~\{\xymatrix@C-5pt{x_1\ar[r] &x_1}\}~,~
\partial x_1^{\perp}~=~\emptyset~.
\end{array}$$
The characteristic map $i:x_0^{\perp} \to X'$ is surjective but not injective, and
$$H_n(x_0^{\perp},\partial x_0^{\perp})~=~
H_n(i(x_0^{\perp}),i(\partial x_0^{\perp}))~=~
\begin{cases} \Z &\hbox{\rm if}~n=1 \\
0&\hbox{\rm if}~n\neq 1~.
\end{cases}$$
\hfill\qed
\end{example}

\begin{example} \label{dual cell3}
Let $X$ be the contractible $\Delta$-set with
one 0-simplex $x_0$, one 1-simplex $x_1$ and one 2-simplex $x_2$,
with non-identity morphisms
$$\xymatrix@C-5pt{x_0 \ar@<0.75ex>[r]  \ar@<-0.75ex>[r] & x_1}~,~
\xymatrix@C-5pt{x_1 \ar@<1.5ex>[r]  \ar[r]  \ar@<-1.5ex>[r]  & x_2}~,~
\xymatrix@C-5pt{x_0 \ar@<1.5ex>[r]  \ar[r]  \ar@<-1.5ex>[r]  & x_2}$$
$$\xymatrix{&x_0 \ar[ddl]_-{\di{x_1}} \ar[ddr]^{\di{x_1}} & \\
& x_2 &\\
x_0 \ar[rr]^{\di{x_1}} && x_0}$$
The realization $\Vert X \Vert$ is the dunce hat (Zeeman \cite{z}).
The barycentric subdivision $X'$ is the $\Delta$-set with three 0-simplices,
eight 1-simplices and six 2-simplices:
$$\begin{array}{l}
{X'}^{(0)}~=~\{x_0,x_1,x_2\}~, \\[1ex]
{X'}^{(1)}~=~\{\xymatrix@C-5pt{x_0 \ar@<0.75ex>[r]  \ar@<-0.75ex>[r] & x_1}\}
\cup
\{\xymatrix@C-5pt{x_1 \ar@<1.5ex>[r]  \ar[r]  \ar@<-1.5ex>[r]  & x_2}\}
\cup
\{\xymatrix@C-5pt{x_0 \ar@<1.5ex>[r]  \ar[r]  \ar@<-1.5ex>[r]  & x_2}\}\\[1ex]
{X'}^{(2)}~=~
\{\xymatrix@C-5pt{x_0 \ar@<0.75ex>[r]  \ar@<-0.75ex>[r] & x_1
\ar@<1.5ex>[r]  \ar[r]  \ar@<-1.5ex>[r]  & x_2}\}
\end{array}$$
$$\xymatrix{&&x_0 \ar@{-}[dd] \ar[dl] \ar[dr] && \\
& x_1 \ar@{->>}[ddl] \ar@{-}[dr] && x_1 \ar@{-}[dl] \ar@{->>}[ddr] \\
&& x_2 \ar@{-}[d] \ar@{-}[dll] \ar@{-}[drr] && \\
x_0 \ar[rr] && x_1 \ar@{->>}[rr] && x_0}$$
The duals and their boundaries are given by
$$\begin{array}{l}
x_0^{\perp}~=~\{\xymatrix@C-5pt{x_0\ar[r] &x_0},~
\xymatrix@C-5pt{x_0 \ar@<0.75ex>[r]  \ar@<-0.75ex>[r] & x_1}~,
\xymatrix@C-5pt{x_0 \ar@<1.5ex>[r]  \ar[r]  \ar@<-1.5ex>[r]  & x_2}\}\\[2ex]
\hskip50pt
\cup\{\xymatrix@C-5pt{x_0 \ar[r] &x_0 \ar@<0.75ex>[r]  \ar@<-0.75ex>[r] & x_1},~
\xymatrix@C-5pt{x_0 \ar[r] & x_0 \ar@<1.5ex>[r]  \ar[r]  \ar@<-1.5ex>[r]  & x_2},~
\xymatrix@C-5pt{x_0 \ar@<0.75ex>[r] \ar@<-0.75ex>[r]
& x_1 \ar@<1.5ex>[r]  \ar[r]  \ar@<-1.5ex>[r]  & x_2}\}\\[2ex]
\hskip150pt \cup\{\xymatrix@C-5pt{x_0 \ar[r] &x_0 \ar@<0.75ex>[r] \ar@<-0.75ex>[r]
& x_1 \ar@<1.5ex>[r]  \ar[r]  \ar@<-1.5ex>[r]  & x_2}\}~,\\[2ex]
\partial x_0^{\perp}~=~\{
\xymatrix@C-5pt{x_0 \ar@<0.75ex>[r]  \ar@<-0.75ex>[r] & x_1}~,
\xymatrix@C-5pt{x_0 \ar@<1.5ex>[r]  \ar[r]  \ar@<-1.5ex>[r]  & x_2}\}
\cup\{\xymatrix@C-5pt{x_0 \ar@<0.75ex>[r] \ar@<-0.75ex>[r]
& x_1 \ar@<1.5ex>[r]  \ar[r]  \ar@<-1.5ex>[r]  & x_2}\}~,~
\Vert \partial x^{\perp}_0 \Vert~\simeq~ S^1 \vee S^1~,\\[2ex]
x_1^{\perp}~=~\{\xymatrix@C-5pt{x_1\ar[r] &x_1}~,~
\xymatrix@C-5pt{x_1 \ar@<1.5ex>[r]  \ar[r]  \ar@<-1.5ex>[r]  & x_2}\}
\cup\{
\xymatrix@C-5pt{x_1\ar[r] &x_1 \ar@<1.5ex>[r]  \ar[r]  \ar@<-1.5ex>[r]  & x_2}\}
~,\\[2ex]
\partial x_1^{\perp}~=~\{\xymatrix@C-5pt{x_1 \ar@<1.5ex>[r]  \ar[r]  \ar@<-1.5ex>[r]  & x_2}\}~,~
\Vert\partial x_1^{\perp}\Vert~\simeq~\{0,1,2\}~,\\[2ex]
x_2^{\perp}~=~\{\xymatrix@C-5pt{x_2 \ar[r] &x_2}\}~,~\partial x_2^{\perp}~=~\emptyset~.
\end{array}$$
The characteristic map $i:x_0^{\perp} \to X'$ is surjective but not injective,
with
$$\Vert i(x^{\perp}_0) \Vert~\simeq~\{*\}~,~
\Vert i(\partial x^{\perp}_0) \Vert~ \simeq~S^1 \vee S^1 ~$$
and
$$H_n(x_0^{\perp},\partial x_0^{\perp})~=~
H_n(i(x_0^{\perp}),i(\partial x_0^{\perp}))~=~
\begin{cases} \Z \oplus \Z&\hbox{\rm if}~n=2 \\
0&\hbox{\rm if}~n\neq 2~.
\end{cases}$$
The characteristic map $i:x_1^{\perp} \to X'$ is neither surjective
nor injective, with
$$\Vert i(x^{\perp}_1) \Vert~\simeq~S^1 \vee S^1~,~
\Vert i(\partial x^{\perp}_1) \Vert~ \simeq ~\{*\}$$
and
$$H_n(x_1^{\perp},\partial x_1^{\perp})~=~
H_n(i(x_1^{\perp}),i(\partial x_1^{\perp}))~=~
\begin{cases} \Z \oplus \Z&\hbox{\rm if}~n=1 \\
0&\hbox{\rm if}~n\neq 1~.
\end{cases}$$
\hfill\qed
\end{example}

\begin{definition}
Given a ring $R$ let $\Mod(R)$ be the additive category of
left $R$-modules.  For $R=\Z$ write $\Mod(\Z)=\Ab$, as usual.
\hfill\qed
\end{definition}

\begin{definition}\label{Delta-chain}
{\rm (Ranicki and Weiss \cite[1.9]{ranickiweiss} for simplicial complexes)}\\
{\rm (i)} The $R$-coefficient simplicial chain complex $\Delta(X';R)$ of the
barycentric subdivision $X'$ of a finite $\Delta$-set $X$ is the
chain complex in $\Mod(R)_*(X)$ with
$$\Delta(X';R)(x)~=~\Delta(x^{\perp},\partial x^{\perp};R)~~,~~\Delta(X';R)[x]~=~\Delta(x^{\perp};R).$$
Compare example~\ref{expl-roundADelta} case (iii). \newline
{\rm (ii)} Let $f:Y \to X'$ be a $\Delta$-map from a
finite $\Delta$-set $Y$ to the barycentric subdivision $X'$ of a $\Delta$-set $X$.
The $R$-coefficient simplicial chain complex $\Delta(Y;R)$ is
the chain complex in $\Mod(R)_*(X)$ with
$$\Delta(Y;R)(x)~=~\Delta(x/f,\partial(x/f);R)~~,~~
\Delta(Y;R)[x]~=~\Delta(x/f;R)~~(x \in X)$$
with $x/f$, $\partial (x/f)$ the $\Delta$-sets defined to fit into
strict pullback squares of $\Delta$-sets
$$\xymatrix{\partial(x/f) \ar[d] \ar[r] &
x/f \ar[d] \ar[r] & Y\ar[d]^{f} \\
\partial x^{\perp} \ar[r] & x^{\perp} \ar[r]^-{i} & X'
}$$
\end{definition}
\hfill\qed

\section{The total complex \label{total complex}}
For a finite chain complex $C$ in $\A_*[X]$,
there is defined a chain complex in $\A^*(X)$, called
the \emph{total complex} of $C$.

\begin{definition} \label{total down}
The \emph{total complex} $\text{Tot}_*C$ of a finite chain complex $C$ in
$\A_*[X]$ is the finite chain complex in $\A^*(X)$ given by
\[
(\text{Tot}_*C)(x)_n=C[x]_{n-|x|}
\]
with differential $d=d_{C[x]}+\sum^{|x|}_{i=0}(-)^{i+|x|}C(\partial_ix\to x)$.
The construction is natural, defining a covariant functor
$$\B(\A)_*[X]\to \B(\A)^*(X)~;~C\mapsto \text{Tot}_*C~.$$
\hfill\qed
\end{definition}

\begin{remark} There is a forgetful functor
$\B(\A)^*_f(X)\to \B(\A)$ taking $C$ in $\B(\A)^*(X)$ to
\[  C(X)= \bigoplus_{x\in X} C(x)~. \]
Compare example~\ref{expl-roundADelta}.
The chain complex $(\text{Tot}_*C)(X)$
in $\A$ is the `realization'
$$\begin{array}{ll}
\bigg(\sum\limits_{x \in X}\Delta(\Delta^{\vert x \vert}) \otimes_{\mathbb{Z}}
C[x]\bigg)/\sim
\end{array}$$
with $\sim$ the equivalence relation generated by
$a \otimes \lambda^*b \sim \lambda_*a \otimes b$ for a morphism $\lambda:y \to z$
in $X$, with $a\in \Delta(\Delta^{\vert y \vert})$,~
$b\in C[z]$. \hfill\qed
\end{remark}

\begin{example} \label{D}
The simplicial chain complex $\Delta(X)$ of a finite $\Delta$-set $X$
is $(\text{Tot}_*C)(X)$ for the chain complex $C$ in
$\Ab_*[X]$ defined by $C[x]~=~\Z$ for all $x$ (a constant functor).
\hfill\qed
\end{example}

\begin{remark} There are evident forgetful functors
\begin{align*}
&\B(\A)_*(X)\to \B(\A)~;~C\mapsto C(X)~,\\
&\B(\A)^*_f(X)\to \B(\A)~;~C\mapsto C(X)~.
\end{align*}
The diagram
$$\xymatrix{
\B(\A)_*(X) \ar[r] \ar[dr] & \B(\A)_*[X]
\ar[r]^{\text{\rm Tot}_*} & \B(\A)_f^*(X) \ar[dl] \\
& \B(\A) & }$$
commutes up to natural chain homotopy equivalence:
for any finite chain complex $C$ in $\A_*(X)$
$$(\text{Tot}_*C)(X)_n~=~\sum\limits_{x \in X} \sum\limits_{x \to y}C(y)_{n-\vert x \vert}~
=~\sum\limits_{y \in X}(\Delta(X/y)\otimes_{\Z}C(y))_n$$
with $X/y$ the $\Delta$-set defined in  \ref{over}, which is contractible.
\hfill\qed
\end{remark}

\begin{proposition}
{\rm (i)} For any objects $M,N$ in $\A_*(X)$ the abelian
group\linebreak $\Hom_{\A_*(X)}(M,N)$ is naturally an object in
$\Ab^*_f(X)$, with
$$\begin{array}{ll}
\Hom_{\A_*(X)}(M,N)(x)&=~\Hom_{\A}(M(x),[N][x])\\[1ex]
&=~\sum\limits_{x \to y} \Hom_{\A}(M(x),N(y))~~(x \in X)~.
\end{array}$$
If $f:M'\to M$, $g:N \to N'$ are morphisms
in $\A_*(X)$ there is induced a morphism in $\Ab^*(X)$
$$\Hom_{\A_*(X)}(M,N) \to \Hom_{\A_*(X)}(M',N')~;~
h \mapsto ghf~.$$
{\rm (ii)}  For any objects $M,N$ in $\A^*_f(X)$ the abelian
group $\Hom_{\A^*(X)}(M,N)$ is naturally
an object in $\Ab_*(X)$, with
$$\begin{array}{ll}
\Hom_{\A^*(X)}(M,N)(x)&=~\Hom_{\A}(M(x),[N][x])\\[1ex]
&=~\sum\limits_{y \to x} \Hom_{\A}(M(x),N(y))~~(x \in X)~.
\end{array}$$
Naturality as in {\rm (i)}.
\end{proposition}
\begin{proof} Immediate from \ref{adjoint}.
\end{proof}

\begin{example} \label{hom}
(i) For a chain complex $C$ in $\Ab_*(X)$ the total complex in
$\Ab^*(X)$ of the corresponding chain complex $[C]$ in $\Ab_*[X]$ is
given by
$$[C]_*[X]~=~\Hom_{\Ab_*(X)}(\Delta(X)^{-*},C)~.$$
(ii) For a chain complex $D$ in $\Ab^*(X)$ the total complex in
$\Ab_*(X)$ of the corresponding chain complex $[D]$ in
$\Ab^*[X]$ is given by
$$[D]^*[X]~=~\Hom_{\Ab^*(X)}(\Delta(X),D)~.$$
\hfill\qed
\end{example}


\section{Chain duality in $L$-theory\label{Duality}}
In general, it is not possible to extend an involution $T:\A
\to \A$ on an additive category $\A$ to the functor
category $\A_*(X)$ for an arbitrary category $X$.  An object in
$\A_*(X)$ is an induced contravariant functor $F:X \to
\A$ and the composite of the contravariant functors
$$\xymatrix{X \ar[r]^-{\di{F}} & \A \ar[r]^-{\di{T}} & \A}$$
is a covariant functor, not a contravariant functor, let alone an induced contravariant functor.
A `chain duality' on $\A$ is essentially an involution on the derived
category of finite chain complexes and chain homotopy classes of chain maps;
an involution on $\A$ is an example of a chain duality.
Given a chain duality on $\A$ we shall now
define a chain duality on the induced functor category $\A_*(X)$, for
any $\Delta$-set $X$, essentially in the same way as was carried out
for a simplicial complex $X$ in \cite{ranicki1}.

\begin{definition} (Ranicki \cite[1.1]{ranicki1}) A {\it chain duality}
$(T,e)$ on an additive category $\A$ is a contravariant additive functor
$$T~:~ \A~ \to ~\B\,(\A)$$
together with a natural transformation
$$e~:~T^2\to 1~: ~\A\to \B\,(\A)$$
such that for each object $M$ in $\A$
\begin{itemize}
\item[(i)] $e(T(M))\circ T(e(M)) = 1 :  T(M) \to T^3(M) \to T(M)~,$
\item[(ii)] $e(M): T^2(M)\to M$ is a chain equivalence.
\end{itemize}
\hfill\qed
\end{definition}
A chain duality $(T,e)$ on $\A$ extends to a contravariant functor
on the bounded chain complex category
$$T~:~\B(\A)~ \to~ \B(\A)~;~C\mapsto T(C)~,$$
using the double complex construction with
$$T(C)_n~=~\sum\limits_{p+q=n} T(C_{-p})_q~~,~~
d_{T(C)}~=~d_{T(C_{-p})} + (-)^qT(d:C_{-p+1}\to C_{-p})~,$$
and $e(C):T^2(C) \to C$ a chain equivalence.
For any objects $M,N$ in an additive category $\A$ there is defined
a $\Z$-module $\Hom_{\A}(M,N)$. Thus for any chain complexes $C,D$
in $\A$ there is defined a $\Z$-module chain complex $\Hom_{\A}(C,D)$,
with
$$\Hom_{\A}(C,D)_n~=~\sum\limits_{q-p=n}\Hom_{\A}(C_p,D_q)~,~
d_{\Hom_{\A}(C,D)}(f)~=~d_D f + (-)^qf d_C~.$$
If $(T,e)$ is a chain duality on $\A$ there is defined a $\Z$-module
chain map
$$\Hom_{\A}(TC,D) \to \Hom_{\A}(TD,C)~;~f \mapsto e(C)T(f)$$
which is a chain equivalence for finite $C$.

\begin{example} An involution $(T,e)$ on $\A$ is a
contravariant functor $T:\A\to \A$
with a natural equivalence $e:T^2\to 1$ such that
for each object $M$ in $\A$
$$e(T(M))~=~T(e(M)^{-1})~:~T^3(M) \to T(M)~.$$
This is essentially the same as a chain duality $(T,e)$ such that
$T(M)$ is a 0-dimensional chain complex for each object $M$ in $\A$.\hfill\qed
\end{example}

\begin{definition}
A {\it chain product} $(\otimes_{\A},b)$ on an additive category
$\A$ is a natural pairing
$$\otimes_{\A}~:~\Ob(\A) \times \Ob(\A) \to
\{\text{$\Z$-module chain complexes}\}~;~
(M,N) \mapsto M\otimes_{\A}N$$
together with a natural chain equivalence
$$b(M,N)~:~M\otimes_{\A}N \to N\otimes_{\A}M$$
such that up to natural isomorphism
$$\begin{array}{l}
(M \oplus M') \otimes_{\A} N~=~(M \otimes_{\A}N) \oplus
(M'\otimes_{\A}N) ~,\\[1ex]
M\otimes_{\A} (N \oplus N') ~=~(M \otimes_{\A}N) \oplus
(M\otimes_{\A}N')
\end{array}$$
and
$$b(N,M) \circ b(M,N)~\simeq~1~:~M \otimes_{\A}N \to M\otimes_{\A}N~.$$
\hfill\qed
\end{definition}

\begin{remark} {\rm The notion of chain product is a linear version
of an `SW-product' in the sense of Weiss and Williams \cite{ww}, where
SW = Spanier-Whitehead.}\\
\hfill\qed
\end{remark}

Given an additive category $\A$ with a chain product $(\otimes_{\A},b)$
and chain complexes $C,D$
in $\A$ let $C \otimes_{\A} D$ be the $\Z$-module chain complex defined by
$$(C\otimes_{\A}D)_n~=~\sum\limits_{p+q+r=n}(C_p \otimes_{\A}D_q)_r~,$$
$$d_{C\otimes_{\A}D}~=~d_{C_p \otimes_{\A} C_q}+(-)^r(1 \otimes_{\A} d_D + (-)^qd_C \otimes_{\A}  1)~.$$
By the naturality of $b$ there is defined a natural chain equivalence
$$b(C,D)~:~C\otimes_{\A}D \to D\otimes_{\A}C~.$$
\begin{proposition} Let $\A$ be an additive category.\\
{\rm (i)} A chain duality $(T,e)$ on $\A$ determines a chain product
$(\otimes_{\A},b)$ on $\A$ by
$$\begin{array}{l}
M \otimes_{\A} N~=~\Hom_{\A}(TM,N)~,\\[1ex]
b(M,N)~:~M \otimes_{\A} N \to N \otimes_{\A}M~;\\[1ex]
\hskip50pt
(f:TM \to N) \mapsto (e(M) \circ T(f): TN \to T^2M \to M)~.
\end{array}$$
{\rm (ii)} If $(\otimes_{\A},b)$ is a chain product on $\A$ such that
$$M \otimes_{\A} N~=~\Hom_{\A}(TM,N)~,~b(M,N)(f)~=~e(M) \circ T(f)$$
for some contravariant additive functor $T:\A \to \B(\A)$ and
natural transformation $e:T^2\to 1:\A\to \B\,(\A)$, then $(T,e)$
is a chain duality on $\A$.
\end{proposition}
\begin{proof} Immediate from the definitions.
\end{proof}

\begin{example} Let $R$ be a ring with an involution
$R \to R;r \mapsto \overline{r}$. Regard a (left) $R$-module $M$
as a right $R$-module by
$$M \times R \to M~;~(x,r) \mapsto \overline{r}x~.$$
Thus for any $R$-modules $M,N$ there is defined a $\Z$-module
$$M\otimes_RN~=~
(M\otimes_{\Z}N)/\{\overline{r}x \otimes y - x\otimes ry\,\vert\,
x \in M, y\in N, r \in R\}$$
with a natural isomorphism
$$b(M,N)~:~M\otimes_RN \to N \otimes_RM~;~x \otimes y \mapsto y \otimes x$$
defining a (0-dimensional) chain product $(\otimes_R,b)$ on the
$R$-module category  $\Mod(R)$. As in \ref{chain} use the
involution on $R$ to define the contravariant duality functor
$$T~:~\Mod(R) \to \Mod(R)~;~M \mapsto M^*=\Hom_R(M,R)$$
with
$$R \times M^* \to M^*~;~(r,f) \mapsto (x \mapsto f(x) \overline{r})~.$$
The natural $\Z$-module morphism defined for any $R$-modules $M,N$ by
$$M \otimes_RN \to \Hom_R(M^*,N)~;~x \otimes y \mapsto (f \mapsto f(x)y)$$
is an isomorphism for f.g. projective $M$.
The $R$-module morphism defined for any $R$-module $M$ by
$$e'(M)~:~M \to M^{**}~;~ x \mapsto (f \mapsto f(x))$$
is an isomorphism for f.g. projective $M$.
Let $\Proj(R) \subset \Mod(R)$ be the full subcategory
of f.g. projective $R$-modules. The natural isomorphisms
$$e(M)~=~e'(M)^{-1}~:~M^{**} \to M$$
define an involution $(T,e)$ on $\Proj(R)$, corresponding to the
restriction to $\Proj(R)$ of the chain product $(\otimes_R,b)$
on $\Mod(R)$.\hfill\qed
\end{example}

\begin{proposition} \label{chain duality on A}
{\rm (Ranicki \cite[5.1,5.9,7]{ranicki1}, Weiss \cite[1.5]{weiss})} \\
A chain duality $(T_{\A},e_{\A})$ on an additive category $\A$
extends to a chain duality $(T_{\A_*(X)},e_{\A_*(X)})$ on $\A_*(X)$,
for any $\Delta$-set $X$
$$\xymatrix{T_{\A_*(X)}~:~\A_*(X) \ar[r] & \A_*[X]
\ar[r]^-{\di{{\rm Tot}_*}}& \B(\A)^*_f(X) \ar[r]^-{\di{T_{\A}}} & \B(\A)_*(X)}$$
where $T_{\A}:\B(\A)^*_f(X) \to \B(\A)_*(X)$
is the extension of the contravariant functor
$$T_{\A}~:~\A^*_f(X)\to
\B(\A)_*(X)~;~M~=~\sum\limits_{x \in X}M(x)\mapsto
T_{\A}(M)~=~\sum\limits_{x \in X}T_{\A}(M(x))~.$$
\end{proposition}

More explicitly, the chain dual of a finite chain complex $C$ in
$\A_*(X)$ is given by
$$T_{\A_*(X)}(C)~=~T_{\A}(\text{Tot}_*C)~,$$
so that
\begin{align*}
T_{\A_*(X)}(C)(x)~&=~T_{\A}(C[x]_{*-\vert x \vert})\\
&=~\sum\limits_{x \to y}T_{\A}(C(y)_{*-\vert x \vert})~~(x \in X)~.
\end{align*}

\begin{example} Let $\A=\A(\mathbb{Z})$, the additive
category of f.g. free abelian groups.\\
(i) For any finite chain complex $C$ in $\A_*(X)$, which we also
view as a (degreewise) induced chain complex $C$ in $\A_*(X)$,
the total complex $\text{Tot}_*(C)$ is given by \ref{hom} to be
$$\Hom_{\A_*(X)}(\Delta(X)^{-*},C)~,$$
so that the chain dual of $C$ is given by
$$T_{\A_*(X)}(C)~=~
\Hom_{\A}(\Hom_{\A_*(X)}(\Delta(X)^{-*},C),\mathbb{Z})~.$$
(ii) As in \ref{Delta-chain} regard the
simplicial chain complex $\Delta(X')$ of the barycentric subdivision $X'$
of a finite $\Delta$-set $X$ as a chain complex in $\A_*(X)$ or in $\A_*[X]$
with
$$\Delta(X')(x)~=~\Delta(x^{\perp},\partial x^{\perp})~,~\Delta(X')[x]~=~\Delta(x^{\perp}) $$
for $x \in X$. The chain dual $T(\Delta(X'))$ is the chain complex in
$\A_*(X)$ with
$$T(\Delta(X'))(x)~=~\Delta(x^{\perp})^{\vert x\vert-*}~~(x \in X)~.$$
\hfill\qed
\end{example}

\begin{remark}
See Fimmel \cite{fimmel} and Woolf \cite{woolf} for Verdier duality
for local coefficient systems on simplicial sets and simplicial complexes.
In particular, \cite{woolf} relates the chain duality of \cite[Chapter 5]{ranicki1}
defined on $\Proj(R)_*(X)$ for a simplicial complex $X$ to the Verdier
duality for sheaves of $R$-module chain complexes over the polyhedron
$\Vert X \Vert$.\hfill\qed
\end{remark}

For any additive category with chain duality ${\mathbb A}$ let
${\mathbb L}_{\bullet}({\mathbb A})$ be the quadratic $L$-theory
$\Omega$-spectrum
defined in Ranicki \cite{ranicki1}, with homotopy groups
$$\pi_n({\mathbb L}_{\bullet}({\mathbb A}))~=~L_n({\mathbb A})~.$$

It was shown in \cite[Chapter 13]{ranicki1} that the covariant functor
$$\{{\rm simplicial~ complexes}\} \to \{\Omega{\rm -spectra}\}~;~
X \mapsto {\mathbb L}_{\bullet}({\mathbb B}({\mathbb A}_*(X)))$$
is an unreduced homology theory, i.e. a covariant functor which
is homotopy invariant, excisive and sends arbitrary disjoint unions to wedges.
More generally~:
\begin{proposition} \label{quad}
{\rm (\cite[13.7]{ranicki1} for simplicial complexes)} \\
{\rm (i)} If $\A$ is an additive category with chain duality and
$X$ is a $\Delta$-set then $\A_*(X)$ is an additive category
with chain duality. \\
{\rm (ii)} The functor
$$\{\Delta{\rm -sets}\} \to \{\Omega{\rm -spectra}\}~;~
X \mapsto L_*(\A,X)~=~\LL_{\bullet}(\B(\A_*(X)))$$
is an unreduced homology theory, that is $L_*(\A,X)=H_*(X;\LL_{\bullet}(A)$.\\
{\rm (iii)} Let $R$ be a ring with involution, so that $\A=\Proj(R)$ is an
additive category of f.g. projective $R$-modules with the duality involution.
If $X$ is a $\Delta$-set and $p:\widetilde{X} \to X$ is a regular cover
with group of covering translations $\pi$ (e.g. the universal cover
with $\pi=\pi_1(X)$) the assembly functor
$$\begin{array}{c}
A~:~\B(R)_*(X) \to \B(R[\pi])~;~ C\mapsto C(\widetilde{X})\\[1ex]
(C(\widetilde{X})~=~\sum\limits_{x\in \widetilde{X}}C(p(x)))
\end{array}$$
is a functor of additive categories with chain duality. The assembly
maps $A$ induced in the $L$-groups fit into an exact sequence
$$\xymatrix@C-10pt{
\dots \ar[r]&H_n(X;\LL_{\bullet}(R)) \ar[r]^-{\di{A}} &
L_n(R[\pi_1(X)]) \ar[r]&\st_n(R,X) \ar[r]&H_{n-1}(X;\LL_{\bullet}(R)) \ar[r]&\dots}$$
with $\st_n(R,X)$ the cobordism group of the $R[\pi_1(X)]$-contractible
$(n-1)$-dimensional quadratic Poincar\'e complexes in $\A_*(X)$.
\end{proposition}
\begin{proof} Exactly as for the simplicial complex case, but using the
$\Delta$-set duals instead of the dual cells!
\end{proof}

\begin{example} {\rm Let $X=S^1$ be the $\Delta$-set of the circle
(\ref{circle1}, \ref{circle2}) with one 0-simplex and one 1-simplex.
Given a ring with involution $R$ let the Laurent polynomial extension ring
$R[z,z^{-1}]$ have the involution $\overline{z}=z^{-1}$.
An $n$-dimensional quadratic Poincar\'e complex in $\Proj(R)_*(S^1)$ is
an $n$-dimensional fundamental quadratic Poincar\'e cobordism over $R$,
with assembly the union $n$-dimensional quadratic Poincar\'e complex over
$R[z,z^{-1}]$, and the assembly maps
$$A~:~H_n(S^1;\LL_{\bullet}(R))~=~L_n(R)\oplus L_{n-1}(R) \to
L_n(R[z,z^{-1}])$$
are isomorphisms modulo the usual $K$-theoretic decorations
(Ranicki \cite[Chapter 24]{ranicki2}.}\hfill\qed
\end{example}

\begin{remark} {\rm Proposition \ref{quad} has an evident analogue for
the symmetric $L$-groups $L^*$.\hfill\qed}
\end{remark}

\providecommand{\bysame}{\leavevmode\hbox to3em{\hrulefill}\thinspace}

\bigskip
\noindent Andrew Ranicki\\
    School of Mathematics\\University of
Edinburgh\\Edinburgh EH9 3JZ\\Scotland,
UK\\
   E-mail: a.ranicki@ed.ac.uk\\

\medskip

\noindent Michael Weiss\\
    School of Mathematical Sciences\\University of
Aberdeen \\Aberdeen AB24 3UE\\ Scotland,
UK \\
   E-mail:  m.weiss@abdn.ac.uk

\end{document}